# SIMULTANEOUS ADAPTATION TO THE MARGIN AND TO COMPLEXITY IN CLASSIFICATION

By Guillaume Lecué

*Université Paris VI*

We consider the problem of adaptation to the margin and to complexity in binary classification. We suggest an exponential weighting aggregation scheme. We use this aggregation procedure to construct classifiers which adapt automatically to margin and complexity. Two main examples are worked out in which adaptivity is achieved in frameworks proposed by Steinwart and Scovel [*Learning Theory. Lecture Notes in Comput. Sci.* **3559** (2005) 279–294. Springer, Berlin; *Ann. Statist.* **35** (2007) 575–607] and Tsybakov [*Ann. Statist.* **32** (2004) 135–166]. Adaptive schemes, like ERM or penalized ERM, usually involve a minimization step. This is not the case for our procedure.

**1. Introduction.** Let $(\mathcal{X}, \mathcal{A})$ be a measurable space. Denote by $D_n$ a sample $((X_i, Y_i))_{i=1,\ldots,n}$ of i.i.d. random pairs of observations where $X_i \in \mathcal{X}$ and $Y_i \in \{-1, 1\}$. Denote by $\pi$ the joint distribution of $(X_i, Y_i)$ on $\mathcal{X} \times \{-1, 1\}$, and $P^X$ the marginal distribution of $X_i$. Let $(X, Y)$ be a random pair distributed according to $\pi$ and independent of the data, and let the component $X$ of the pair be observed. The problem of statistical learning in classification (pattern recognition) consists of predicting the corresponding value $Y \in \{-1, 1\}$.

A *prediction rule* is a measurable function $f : \mathcal{X} \longmapsto \{-1, 1\}$. The *misclassification error* associated with $f$ is

$$R(f) = \mathbb{P}(Y \neq f(X)).$$

It is well known (see, e.g., Devroye, Györfi and Lugosi [15]) that

$$\min_f R(f) = R(f^*) = R^*, \qquad \text{where } f^*(x) = \text{sign}(2\eta(x) - 1)$$

and $\eta$ is the *a posteriori probability* defined by

$$\eta(x) = \mathbb{P}(Y = 1 | X = x),$$

for all $x \in \mathcal{X}$ [where $\text{sign}(y)$ denotes the sign of $y \in \mathbb{R}$ with the convention $\text{sign}(0) = 1$]. The prediction rule $f^*$ is called the *Bayes rule* and $R^*$ is called the *Bayes risk*. A *classifier* is a function, $\hat{f}_n = \hat{f}_n(X, D_n)$, measurable with respect to $D_n$ and $X$ with values in $\{-1, 1\}$ that assigns to every sample $D_n$ a prediction rule $\hat{f}_n(\cdot, D_n) : \mathcal{X} \longmapsto \{-1, 1\}$. A key characteristic of $\hat{f}_n$ is the *generalization error* $\mathbb{E}[R(\hat{f}_n)]$, where

$$R(\hat{f}_n) = \mathbb{P}(Y \neq \hat{f}_n(X) | D_n).$$

The aim of statistical learning is to construct a classifier $\hat{f}_n$ such that $\mathbb{E}[R(\hat{f}_n)]$ is as close to $R^*$ as possible. Accuracy of a classifier $\hat{f}_n$ is measured by the value $\mathbb{E}[R(\hat{f}_n)] - R^*$, called the *excess risk* of $\hat{f}_n$.

The classical approach due to Vapnik and Chervonenkis (see, e.g., [15]) consists of searching for a classifier that minimizes the *empirical risk*

$$(1.1) \qquad R_n(f) = \frac{1}{n} \sum_{i=1}^{n} \mathbb{1}_{(Y_i f(X_i) \leq 0)},$$

over all prediction rules $f$ in a source class $\mathcal{F}$, where $\mathbb{1}_A$ denotes the indicator of the set $A$. Minimizing the empirical risk (1.1) is computationally intractable for many sets $\mathcal{F}$ of classifiers, because this functional is neither convex nor continuous. Nevertheless, we might base a tractable estimation procedure on minimization of a convex surrogate $\phi$ for the loss (Cortes and Vapnik [13], Freund and Schapire [17], Lugosi and Vayatis [28], Friedman, Hastie and Tibshirani [18] and Bühlmann and Yu [7]). It has recently been shown that these classification methods often give classifiers with small Bayes risk (Blanchard, Lugosi and Vayatis [5] and Steinwart and Scovel [38, 39]). The main idea is that the sign of the minimizer of $A^{(\phi)}(f) = \mathbb{E}[\phi(Y f(X))]$, the $\phi$-*risk*, where $\phi$ is a convex loss function and $f$ a real-valued function, is in many cases equal to the Bayes classifier $f^*$. Therefore, minimizing $A_n^{(\phi)}(f) = \frac{1}{n} \sum_{i=1}^{n} \phi(Y_i f(X_i))$, the *empirical $\phi$-risk*, and taking $\hat{f}_n = \text{sign}(\hat{F}_n)$ where $\hat{F}_n \in \text{Arg min}_{f \in \mathcal{F}} A_n^{(\phi)}(f)$, leads to an approximation for $f^*$. Here, $\text{Arg min}_{f \in \mathcal{F}} P(f)$, for a functional $P$, denotes the set of all $f \in \mathcal{F}$ such that $P(f) = \min_{f \in \mathcal{F}} P(f)$. Schapire, Freund, Bartlett and Lee [36], Lugosi and Vayatis [28], Blanchard, Lugosi and Vayatis [5], Zhang [48], Steinwart and Scovel [38, 39] and Bartlett, Jordan and McAuliffe [2] give results on statistical properties of classifiers obtained by minimization of such a convex risk. A wide variety of classification methods in machine learning are based on this idea, in particular, on using the convex loss associated with support vector machines (Cortes and Vapnik [13] and Schölkopf and Smola [37]),

$$\phi(x) = (1-x)_+,$$

called the *hinge-loss*, where $z_+ = \max(0, z)$ denotes the positive part of $z \in \mathbb{R}$. Denote by

$$A(f) = \mathbb{E}[(1 - Yf(X))_+]$$

the *hinge risk* of $f : \mathcal{X} \longmapsto \mathbb{R}$ and set

(1.2) $$A^* = \inf_f A(f),$$

where the infimum is taken over all measurable functions $f$. We will call $A^*$ the *optimal hinge risk*. One may verify that the Bayes rule $f^*$ attains the infimum in (1.2) and Lin [27] and Zhang [48] have shown that

(1.3) $$R(f) - R^* \leq A(f) - A^*,$$

for all measurable functions $f$ with values in $\mathbb{R}$. Thus, minimization of $A(f) - A^*$, the *excess hinge risk*, provides a reasonable alternative for minimization of excess risk.

The difficulty of classification is closely related to the behavior of the a posteriori probability $\eta$. Mammen and Tsybakov [31], for the problem of discriminant analysis which is close to our classification problem, and Tsybakov [42] have introduced an assumption on the closeness of $\eta$ to $1/2$, called the *margin assumption* (or *low noise assumption*). Under this assumption, the risk of a minimizer of the empirical risk over some fixed class $\mathcal{F}$ converges to the minimum risk over the class with a *fast rate*, namely, faster than $n^{-1/2}$. In fact, with no assumption on the joint distribution $\pi$, the convergence rate of the excess risk is not faster than $n^{-1/2}$ (cf. Devroye, Györfi and Lugosi [15]). However, under the margin assumption, it can be as fast as $n^{-1}$. Minimizing a penalized empirical hinge risk, under this assumption, also leads to fast convergence rates (Blanchard, Bousquet and Massart [4], Steinwart and Scovel [38, 39]). Massart [32], Massart and Nédélec [34] and Massart [33] also obtain results that can lead to fast rates in classification using penalized empirical risk in the special case of a low noise assumption. Audibert and Tsybakov [1] show that fast rates can be achieved for plug-in classifiers.

In this paper we consider the problem of adaptive classification. Mammen and Tsybakov [31] have shown that fast rates depend on both the *margin parameter* $\kappa$ and complexity $\rho$ of the class of candidate sets for $\{x \in \mathcal{X} : \eta(x) \geq 1/2\}$. Their results were nonadaptive, supposing that $\kappa$ and $\rho$ are known. Tsybakov [42] suggested an adaptive classifier that attains fast optimal rates, up to a logarithmic factor, without knowing $\kappa$ and $\rho$. Tsybakov and van de Geer [43] suggest a penalized empirical risk minimization classifier that adaptively attains, up to a logarithmic factor, the same fast optimal rates of convergence. Tarigan and van de Geer [40] extend this result to $l_1$-penalized empirical hinge risk minimization. Koltchinskii [23] uses Rademacher averages to get a similar result without the logarithmic factor. Related work is that of Koltchinskii [22], Koltchinskii and Panchenko [24] and Lugosi and Wegkamp [29].

Note that the existing papers on fast rates either suggest classifiers that can be easily implemented but are nonadaptive, or adaptive schemes that are hard to apply in practice and/or do not achieve the minimax rates (they pay a price for adaptivity). The aim of the present paper is to suggest and to analyze an exponential weighting aggregation scheme which does not require a minimization step, unlike other adaptation schemes such as ERM (Empirical Risk Minimization) and penalized ERM, and which does not pay a price for adaptivity. This scheme is used first to construct minimax adaptive classifiers (cf. Theorem 3.1) and second to construct easily implemented classifiers that are adaptive simultaneously to complexity and to the margin parameters and which achieve the fast rates.

The paper is organized as follows. In Section 2 we prove an oracle inequality which corresponds to the adaptation step of the procedure that we suggest. In Section 3 we apply the oracle inequality to two types of classifiers, one of which is constructed by minimization on sieves (as in Tsybakov [42]), and which gives an adaptive classifier which attains fast optimal rates without a logarithmic factor, and the other which is based on support vector machines (SVM), following Steinwart and Scovel [38, 39]. The latter is realized as a computationally feasible procedure and it adaptively attains fast rates of convergence. In particular, we suggest a method of adaptive choice of the parameter of $L1$-SVM classifiers with Gaussian RBF kernels. Proofs are given in Section 4.

**2. Oracle inequalities.** In this section we give an oracle inequality showing that a specifically defined convex combination of classifiers mimics the best classifier in a given finite set.

Suppose that we have $M \geq 2$ different classifiers $\hat{f}_1, \ldots, \hat{f}_M$ taking values in $\{-1, 1\}$. The problem of model selection type aggregation, as studied in Nemirovski [35], Yang [46, 47], Catoni [11] and Tsybakov [41], consists in construction of a new classifier $\tilde{f}_n$ (called *aggregate*) which is approximatively at least as good, with respect to the excess risk, as the best among $\hat{f}_1, \ldots, \hat{f}_M$. In most of these papers the aggregation is based on a splitting of the sample into two independent subsamples $D_m^1$ and $D_l^2$ of sizes $m$ and $l$, respectively, where $m \gg l$ and $m + l = n$. The first subsample $D_m^1$ is used to construct the classifiers $\hat{f}_1, \ldots, \hat{f}_M$ and the second subsample $D_l^2$ is used to aggregate them, that is, to construct a new classifier that mimics in a certain sense the behavior of the best among the classifiers $\hat{f}_i$.

In this section we will not consider the sample splitting and will concentrate only on the construction of aggregates (following Nemirovski [35], Juditsky and Nemirovski [20], Tsybakov [41], Birgé [3] and Bunea, Tsybakov and Wegkamp [10]). Thus, the first subsample is fixed, and instead of classifiers $\hat{f}_1, \ldots, \hat{f}_M$, we have fixed prediction rules $f_1, \ldots, f_M$. Rather than work with a part of the initial sample, we will suppose, for notational simplicity, that the whole sample $D_n$ of size $n$ is used for the aggregation step instead of a subsample $D_l^2$.

Our procedure uses exponential weights. The idea of exponential weights is well known; see, for example, Buckland, Burnham and Augustin [8], Yang [47], Catoni [11], Hartigan [19] and Leung and Barron [26]. This procedure has been widely used in on-line prediction; see, for example, Vovk [45] and Cesa-Bianchi and Lugosi [12]. We consider the following aggregate which is a convex combination with exponential weights of $M$ classifiers:

$$\tilde{f}_n = \sum_{j=1}^{M} w_j^{(n)} f_j, \tag{2.1}$$

where

$$w_j^{(n)} = \frac{\exp(\sum_{i=1}^{n} Y_i f_j(X_i))}{\sum_{k=1}^{M} \exp(\sum_{i=1}^{n} Y_i f_k(X_i))} \qquad \forall j = 1, \ldots, M. \tag{2.2}$$

Since $f_1, \ldots, f_M$ take their values in $\{-1, 1\}$, we have

$$w_j^{(n)} = \frac{\exp(-nA_n(f_j))}{\sum_{k=1}^{M} \exp(-nA_n(f_k))}, \tag{2.3}$$

for all $j \in \{1, \ldots, M\}$, where

$$A_n(f) = \frac{1}{n} \sum_{i=1}^{n} (1 - Y_i f(X_i))_+ \tag{2.4}$$

is the empirical analog of the hinge risk. Since $A_n(f_j) = 2R_n(f_j)$ for all $j = 1, \ldots, M$, these weights can be written in terms of the empirical risks of the $f_j$'s,

$$w_j^{(n)} = \frac{\exp(-2nR_n(f_j))}{\sum_{k=1}^{M} \exp(-2nR_n(f_k))} \qquad \forall j = 1, \ldots, M.$$

The aggregation procedure defined by (2.1) with weights (2.3) does not need any minimization algorithm in contrast to the ERM procedure. Moreover, the following proposition shows that this exponential weighting aggregation scheme has theoretical properties similar to those of the ERM procedure, up to the residual $(\log M)/n$. In what follows, the aggregation procedure defined by (2.1) with exponential weights (2.3) is called the Aggregation procedure with Exponential Weights and is denoted by AEW.

PROPOSITION 2.1. *Let $M \geq 2$ be an integer and $f_1, \ldots, f_M$ be $M$ prediction rules on $\mathcal{X}$. For any integers $n$, the AEW procedure $\tilde{f}_n$ satisfies*

$$(2.5) \qquad A_n(\tilde{f}_n) \leq \min_{i=1,\ldots,M} A_n(f_i) + \frac{\log(M)}{n}.$$

Obviously, inequality (2.5) is satisfied when $\tilde{f}_n$ is the ERM aggregate defined by

$$\tilde{f}_n \in \operatorname{Arg} \min_{f \in \{f_1,\ldots,f_M\}} R_n(f).$$

It is a convex combination of $f_j$'s with weights $w_j = 1$ for one $j \in \operatorname{Arg\,min}_j R_n(f_j)$ and 0 otherwise.

We will use the following assumption (cf. Mammen and Tsybakov [31] and Tsybakov [42]) that will allow us to get fast learning rates for the classifiers that we aggregate.

ASSUMPTION (MA1) [*Margin (or low noise) assumption*]. The probability distribution $\pi$ on the space $\mathcal{X} \times \{-1,1\}$ satisfies the margin assumption (MA1)$(\kappa)$ with margin parameter $1 \leq \kappa < +\infty$ if there exists $c > 0$ such that

$$(2.6) \qquad \mathbb{E}\{|f(X) - f^*(X)|\} \leq c(R(f) - R^*)^{1/\kappa},$$

for all measurable functions $f$ with values in $\{-1,1\}$.

We first give the following proposition which is valid not necessarily for the particular choice of weights given in (2.2).

PROPOSITION 2.2. *Let Assumption (MA1)$(\kappa)$ hold with some $1 \leq \kappa < +\infty$. Assume that there exist two positive numbers $a \geq 1, b$ such that $M \geq an^b$. Let $w_1, \ldots, w_M$ be $M$ statistics measurable w.r.t. the sample $D_n$, such that $w_j \geq 0$, for all $j = 1, \ldots, M$, and $\sum_{j=1}^M w_j = 1$ ($\pi^{\otimes n}$-a.s.). Define $\tilde{f}_n = \sum_{j=1}^M w_j f_j$, where $f_1, \ldots, f_M$ are prediction rules. There exists a constant $C_0 > 0$ such that*

$$(1 - (\log M)^{-1/4})\mathbb{E}[A(\tilde{f}_n) - A^*]$$
$$\leq \mathbb{E}[A_n(\tilde{f}_n) - A_n(f^*)] + C_0 n^{-\kappa/(2\kappa-1)}(\log M)^{7/4},$$

*where $f^*$ is the Bayes rule. For instance, we can take $C_0 = 10 + ca^{-1/(2b)} + a^{-1/b}\exp[(b(8c/6)^2) \vee (((8c/3) \vee 1)/b)^2]$.*

As a consequence, we obtain the following oracle inequality.

THEOREM 2.3. *Let Assumption* (MA1)($\kappa$) *hold with some* $1 \leq \kappa < +\infty$. *Assume that there exist two positive numbers* $a \geq 1, b$ *such that* $M \geq an^b$. *Let* $\tilde{f}_n$ *satisfy* (2.5), *for instance, the AEW or the ERM procedure. Then* $\tilde{f}_n$ *satisfies*

$$(2.7) \quad \begin{aligned} \mathbb{E}[R(\tilde{f}_n) - R^*] \\ \leq \left(1 + \frac{2}{\log^{1/4}(M)}\right) \left\{ 2 \min_{j=1,\ldots,M} (R(f_j) - R^*) + C_0 \frac{\log^{7/4}(M)}{n^{\kappa/(2\kappa-1)}} \right\} \end{aligned}$$

*for all integers* $n \geq 1$, *where* $C_0 > 0$ *appears in Proposition* 2.2.

REMARK 2.1. The factor 2 multiplying $\min_{j=1,\ldots,M}(R(f_j) - R^*)$ in (2.7) is due to the relation between the hinge excess risk and the usual excess risk [cf. inequality (1.3)]. The hinge-loss is more adapted for our convex aggregate, since we have the same statement without this factor, namely,

$$\mathbb{E}[A(\tilde{f}_n) - A^*] \leq \left(1 + \frac{2}{\log^{1/4}(M)}\right) \left\{ \min_{j=1,\ldots,M}(A(f_j) - A^*) + C_0 \frac{\log^{7/4}(M)}{n^{\kappa/(2\kappa-1)}} \right\}.$$

Moreover, linearity of the hinge-loss on $[-1, 1]$ leads to

$$\min_{j=1,\ldots,M}(A(f_j) - A^*) = \min_{f \in Conv}(A(f) - A^*),$$

where *Conv* is the convex hull of the set $\{f_j : j = 1, \ldots, M\}$. Therefore, the excess hinge risk of $\tilde{f}_n$ is approximately the same as one of the best convex combinations of $f_j$'s.

REMARK 2.2. For a convex loss function $\phi$, consider the empirical $\phi$-risk $A_n^{(\phi)}(f)$. Our proof implies that the aggregate

$$\tilde{f}_n^{(\phi)}(x) = \sum_{j=1}^{M} w_j^{\phi} f_j(x) \qquad \text{with } w_j^{\phi} = \frac{\exp(-nA_n^{(\phi)}(f_j))}{\sum_{k=1}^{M} \exp(-nA_n^{(\phi)}(f_k))}, \ \forall j = 1, \ldots, M,$$

satisfies the inequality (2.5) with $A_n^{(\phi)}$ in place of $A_n$.

We consider next a recursive analog of the aggregate (2.1). It is close to the one suggested by Yang [46] for density aggregation under Kullback loss and by Catoni [11] and Bunea and Nobel [9] for the regression model with squared loss. It can be also viewed as a particular case of the mirror descent algorithm suggested in Juditsky, Nazin, Tsybakov and Vayatis [21]. We consider

$$(2.8) \qquad \bar{f}_n = \frac{1}{n} \sum_{k=1}^{n} \tilde{f}_k = \sum_{j=1}^{M} \bar{w}_j f_j,$$

where

(2.9) $$\bar{w}_j = \frac{1}{n}\sum_{k=1}^n w_j^{(k)} = \frac{1}{n}\sum_{k=1}^n \frac{\exp(-kA_k(f_j))}{\sum_{l=1}^M \exp(-kA_k(f_l))},$$

for all $j = 1, \ldots, M$, where $A_k(f) = (1/k)\sum_{i=1}^k (1 - Y_i f(X_i))_+$ is the empirical hinge risk of $f$ and $w_j^{(k)}$ is the weight defined in (2.2) for the first $k$ observations. This aggregate is especially useful for the on-line framework. The following theorem says that it has the same theoretical properties as the aggregate (2.1).

THEOREM 2.4. *Let Assumption* (MA1)($\kappa$) *hold with some* $1 \leq \kappa < +\infty$. *Assume that there exist two positive numbers* $a \geq 1, b$ *such that* $M \geq an^b$. *Then the convex aggregate* $\bar{f}_n$ *defined by* (2.8) *satisfies*

$$\mathbb{E}[R(\bar{f}_n) - R^*] \leq \left(1 + \frac{2}{\log^{1/4}(M)}\right)\left\{2\min_{j=1,\ldots,M}(R(f_j) - R^*)\right.$$
$$\left. + C_0\gamma(n,\kappa)\log^{7/4}(M)\right\},$$

*for all integers* $n \geq 1$, *where* $C_0 > 0$ *appears in Proposition* 2.2 *and* $\gamma(n,\kappa)$ *is equal to* $((2\kappa - 1)/(\kappa - 1))n^{-\kappa/(2\kappa-1)}$ *if* $\kappa > 1$ *and to* $(\log n)/n$ *if* $\kappa = 1$.

REMARK 2.3. For all $k \in \{1, \ldots, n-1\}$, less observations are used to construct $\tilde{f}_k$ than to construct $\tilde{f}_n$; thus, intuitively, we expect that $\tilde{f}_n$ will learn better than $\tilde{f}_k$. In view of (2.8), $\bar{f}_n$ is an average of aggregates whose performances are, a priori, worse than those of $\tilde{f}_n$; therefore, its expected learning properties are presumably worse than those of $\tilde{f}_n$. An advantage of the aggregate $\bar{f}_n$ is its recursive construction, but the risk behavior of $\tilde{f}_n$ seems to be better than that of $\bar{f}_n$. In fact, it is easy to see that Theorem 2.4 is satisfied for any aggregate $\bar{f}_n = \sum_{k=1}^n w_k \tilde{f}_k$, where $w_k \geq 0$ and $\sum_{k=1}^n w_k = 1$ with $\gamma(n,\kappa) = \sum_{k=1}^n w_k k^{-\kappa/(2\kappa-1)}$, and the remainder term is minimized for $w_j = 1$ when $j = n$ and 0 elsewhere, that is, for $\bar{f}_n = \tilde{f}_n$.

REMARK 2.4. In this section we have dealt only with the aggregation step. But the construction of classifiers has to take place prior to this step. This requires a split of the sample as discussed at the beginning of this section. The main drawback of this method is that only a part of the sample is used for the initial estimation. However, by using different splits of the sample and taking the average of the aggregates associated with each of them, we get a more balanced classifier which does not depend on a particular split. Since the hinge loss is linear on $[-1, 1]$, we have the same result as in Theorem 2.3 and Theorem 2.4 for an average of aggregates of the form (2.1) and (2.8), respectively, for averaging over different splits of the sample.

AGGREGATION OF CLASSIFIERS  9

**3. Adaptation to the margin and to complexity.** In Steinwart and Scovel [38, 39] and Tsybakov [42] two concepts of complexity are used. In this section we show that combining classifiers used by Tsybakov [42] or the $L$1-SVM classifiers of Steinwart and Scovel [38, 39] with our aggregation method leads to classifiers that are adaptive both to the margin parameter and to the complexity in the two cases. Results are established for the first method of aggregation defined in (2.1), but they are also valid for the recursive aggregate defined in (2.8).

We use a sample splitting to construct our aggregate. The first subsample $D_m^1 = ((X_1, Y_1), \ldots, (X_m, Y_m))$, where $m = n - l$ and $l = \lceil an/\log n \rceil$ for a constant $a > 0$, is implemented to construct classifiers and the second subsample $D_l^2 = ((X_{m+1}, Y_{m+1}), \ldots, (X_n, Y_n))$ is implemented to aggregate them by the procedure (2.1).

3.1. *Adaptation in the framework of Tsybakov.* Here we take $\mathcal{X} = \mathbb{R}^d$. Introduce the following pseudo-distance, and its empirical analogue, between the sets $G, G' \subseteq \mathcal{X}$:

$$d_\Delta(G, G') = P^X(G \Delta G'), \qquad d_{\Delta, e}(G, G') = \frac{1}{n} \sum_{i=1}^n \mathbb{1}_{(X_i \in G \Delta G')},$$

where $G \Delta G'$ is the symmetric difference between the sets $G$ and $G'$. If $\mathcal{Y}$ is a class of subsets of $\mathcal{X}$, denote by $\mathcal{H}_B(\mathcal{Y}, \delta, d_\Delta)$ the $\delta$-*entropy with bracketing of* $\mathcal{Y}$ *for the pseudo-distance* $d_\Delta$ (cf. van de Geer [44], page 16). We say that $\mathcal{Y}$ has a *complexity bound* $\rho > 0$ if there exists a constant $A > 0$ such that

$$\mathcal{H}_B(\mathcal{Y}, \delta, d_\Delta) \leq A \delta^{-\rho} \qquad \forall 0 < \delta \leq 1.$$

Various examples of classes $\mathcal{Y}$ having this property can be found in Dudley [16], Korostelëv and Tsybakov [25] and Mammen and Tsybakov [30].

Let $(\mathcal{G}_\rho)_{\rho_{\min} \leq \rho \leq \rho_{\max}}$ be a collection of classes of subsets of $\mathcal{X}$, where $\mathcal{G}_\rho$ has a complexity bound $\rho$, for all $\rho_{\min} \leq \rho \leq \rho_{\max}$. This collection corresponds to a priori knowledge on $\pi$ that the set $G^* = \{x \in \mathcal{X} : \eta(x) > 1/2\}$ lies in one of these classes (typically we have $\mathcal{G}_\rho \subset \mathcal{G}_{\rho'}$ if $\rho \leq \rho'$). The aim of adaptation to the margin and complexity is to propose $\tilde{f}_n$, a classifier free of $\kappa$ and $\rho$ such that, if $\pi$ satisfies (MA1)($\kappa$) and $G^* \in \mathcal{G}_\rho$, then $\tilde{f}_n$ learns with the optimal rate $n^{-\kappa/(2\kappa + \rho - 1)}$ (optimality has been established in Mammen and Tsybakov [31]), and this property holds for all values of $\kappa \geq 1$ and $\rho_{\min} \leq \rho \leq \rho_{\max}$. Following Tsybakov [42], we introduce the following assumption on the collection $(\mathcal{G}_\rho)_{\rho_{\min} \leq \rho \leq \rho_{\max}}$.

ASSUMPTION (A1) (*Complexity assumption*). Assume that $0 < \rho_{\min} < \rho_{\max} < 1$ and the $\mathcal{G}_\rho$'s are classes of subsets of $\mathcal{X}$ such that $\mathcal{G}_\rho \subseteq \mathcal{G}_{\rho'}$ for $\rho_{\min} \leq \rho < \rho' \leq \rho_{\max}$ and the class $\mathcal{G}_\rho$ has complexity bound $\rho$. For any integer $n$, we define $\rho_{n,j} = \rho_{\min} + \frac{j}{N(n)}(\rho_{\max} - \rho_{\min})$, $j = 0, \ldots, N(n)$, where $N(n)$ satisfies $A'_0 n^{b'} \leq N(n) \leq A_0 n^b$, for some finite $b \geq b' > 0$ and $A_0, A'_0 > 0$. Assume that for all $n \in \mathbb{N}$:

 (i) for all $j = 0, \ldots, N(n)$, there exists $\mathcal{N}_n^j$, an $\varepsilon$-net on $\mathcal{G}_{\rho_{n,j}}$ for the pseudo-distance $d_\Delta$ or $d_{\Delta,e}$, where $\varepsilon = a_j n^{-1/(1+\rho_{n,j})}$, $a_j > 0$ and $\max_j a_j < +\infty$;

 (ii) $\mathcal{N}_n^j$ has complexity bound $\rho_{n,j}$, for $j = 0, \ldots, N(n)$.

The first subsample $D_m^1$ is used to construct the ERM classifiers $\hat{f}_m^j(x) = 2\mathbb{1}_{\hat{G}_m^j}(x) - 1$, where $\hat{G}_m^j \in \operatorname{Arg\,min}_{G \in \mathcal{N}_m^j} R_m(2\mathbb{1}_G - 1)$ for all $j = 0, \ldots, N(m)$, and the second subsample $D_l^2$ is used to construct the exponential weights of the aggregation procedure,

$$w_j^{(l)} = \frac{\exp(-lA^{[l]}(\hat{f}_m^j))}{\sum_{k=1}^{N(m)} \exp(-lA^{[l]}(\hat{f}_m^k))} \qquad \forall j = 0, \ldots, N(m),$$

where $A^{[l]}(f) = (1/l)\sum_{i=m+1}^n (1 - Y_i f(X_i))_+$ is the empirical hinge risk of $f : \mathcal{X} \longmapsto \mathbb{R}$ based on the subsample $D_l^2$. We consider

$$(3.1) \qquad \tilde{f}_n(x) = \sum_{j=0}^{N(m)} w_j^{(l)} \hat{f}_m^j(x) \qquad \forall x \in \mathcal{X}.$$

The construction of the $\hat{f}_m^j$'s does not depend on the margin parameter $\kappa$.

THEOREM 3.1. *Let $(\mathcal{G}_\rho)_{\rho_{\min} \leq \rho \leq \rho_{\max}}$ be a collection of classes satisfying Assumption* (A1). *Then the aggregate defined in* (3.1) *satisfies*

$$\sup_{\pi \in \mathcal{P}_{\kappa,\rho}} \mathbb{E}[R(\tilde{f}_n) - R^*] \leq C n^{-\kappa/(2\kappa+\rho-1)} \qquad \forall n \geq 1,$$

*for all $1 \leq \kappa < +\infty$ and all $\rho \in [\rho_{\min}, \rho_{\max}]$, where $C > 0$ is a constant depending only on $a, b, b', A, A_0, A'_0, \rho_{\min}, \rho_{\max}$ and $\kappa$, and $\mathcal{P}_{\kappa,\rho}$ is the set of all probability measures $\pi$ on $\mathcal{X} \times \{-1, 1\}$ such that Assumption* (MA1)$(\kappa)$ *is satisfied and $G^* \in \mathcal{G}_\rho$.*

3.2. *Adaptation in the framework of Steinwart and Scovel.*

3.2.1. *The case of a continuous kernel.* Steinwart and Scovel [38] have obtained fast learning rates for SVM classifiers depending on three parameters, the *margin parameter* $0 \leq \alpha < +\infty$, the complexity exponent $0 < p \leq 2$ and the approximation exponent $0 \leq \beta \leq 1$. The margin assumption was first introduced in Mammen and Tsybakov [31] for the problem of discriminant analysis and in Tsybakov [42] for the classification problem, in the following way:

ASSUMPTION (MA2) [*Margin (or low noise) assumption*]. The probability distribution $\pi$ on the space $\mathcal{X} \times \{-1, 1\}$ satisfies the margin assumption (MA2)$(\alpha)$ with margin parameter $0 \leq \alpha < +\infty$ if there exists $c_0 > 0$ such that

$$(3.2) \qquad \mathbb{P}(|2\eta(X) - 1| \leq t) \leq c_0 t^\alpha \qquad \forall t > 0.$$

As shown in Boucheron, Bousquet and Lugosi [6], the margin Assumptions (MA1)$(\kappa)$ and (MA2)$(\alpha)$ are equivalent with $\kappa = \frac{1+\alpha}{\alpha}$ for $\alpha > 0$.

Let $\mathcal{X}$ be a compact metric space. Let $H$ be a reproducing kernel Hilbert space (RKHS) over $\mathcal{X}$ (see, e.g., Cristianini and Shawe–Taylor [14] and Schölkopf and Smola [37]) and $B_H$ its closed unit ball. Denote by $\mathcal{N}(B_H, \varepsilon, L_2(P_n^X))$ the $\varepsilon$-*covering number* of $B_H$ w.r.t. the canonical distance of $L_2(P_n^X)$, the $L_2$-space w.r.t. the empirical measure, $P_n^X$, on $X_1, \ldots, X_n$. Introduce the following assumptions as in Steinwart and Scovel [38]:

ASSUMPTION (A2). There exist $a_0 > 0$ and $0 < p \leq 2$ such that, for any integer $n$,

$$\sup_{D_n \in (\mathcal{X} \times \{-1,1\})^n} \log \mathcal{N}(B_H, \varepsilon, L_2(P_n^X)) \leq a_0 \varepsilon^{-p} \qquad \forall \varepsilon > 0.$$

Note that the supremum is taken over all samples of size $n$ and the bound is assuming for any $n$. Every RKHS satisfies (A2) with $p = 2$ (cf. Steinwart and Scovel [38]). We define the *approximation error function* of the $L1$-SVM as $a(\lambda) \stackrel{\text{def}}{=} \inf_{f \in H}(\lambda \|f\|_H^2 + A(f)) - A^*$.

ASSUMPTION (A3). The RKHS $H$ approximates $\pi$ with exponent $0 \leq \beta \leq 1$, if there exists a constant $C_0 > 0$ such that $a(\lambda) \leq C_0 \lambda^\beta$, $\forall \lambda > 0$.

Note that every RKHS approximates every probability measure with exponent $\beta = 0$ and the other extremal case $\beta = 1$ is equivalent to the fact that the Bayes classifier $f^*$ belongs to the RKHS (cf. Steinwart and Scovel [38]). Furthermore, $\beta > 1$ only for probability measures such that $\mathbb{P}(\eta(X) = 1/2) = 1$ (cf. Steinwart and Scovel [38]). If (A2) and (A3) hold, the parameter $(p, \beta)$ can be considered as a complexity parameter characterizing $\pi$ and $H$.

Let $H$ be an RKHS with a continuous kernel on $\mathcal{X}$ satisfying (A2) with parameter $0 < p < 2$. Define the $L1$-SVM classifier by

$$(3.3) \qquad \hat{f}_n^\lambda = \text{sign}(\hat{F}_n^\lambda), \text{ where } \hat{F}_n^\lambda \in \text{Arg}\min_{f \in H}(\lambda \|f\|_H^2 + A_n(f));$$

$\lambda > 0$ is called the *regularization parameter*. Assume that the probability measure $\pi$ belongs to the set $\mathcal{Q}_{\alpha,\beta}$ of all probability measures on $\mathcal{X} \times \{-1, 1\}$ satisfying Assumption (MA2)($\alpha$) with $\alpha \geq 0$ and (A3) with complexity parameter $(p,\beta)$, where $0 < \beta \leq 1$. It has been shown in Steinwart and Scovel [38] that the $L1$-SVM classifier, $\hat{f}_n^{\lambda_n^{\alpha,\beta}}$, where the regularization parameter is $\lambda_n^{\alpha,\beta} = n^{-4(\alpha+1)/(2\alpha+p\alpha+4)(1+\beta)}$, satisfies the following excess risk bound: for any $\varepsilon > 0$, there exists $C > 0$ depending only on $\alpha, p, \beta$ and $\varepsilon$ such that

$$(3.4) \quad \mathbb{E}[R(\hat{f}_n^{\lambda_n^{\alpha,\beta}}) - R^*] \leq C n^{-4\beta(\alpha+1)/((2\alpha+p\alpha+4)(1+\beta))+\varepsilon} \qquad \forall n \geq 1.$$

We remark that if $\beta = 1$, that is, $f^* \in H$, then the learning rate in (3.4) is (up to an $\varepsilon$) $n^{-2(\alpha+1)/(2\alpha+p\alpha+4)}$, which is a fast rate since $2(\alpha+1)/(2\alpha+p\alpha+4) \in [1/2, 1)$.

To construct the classifier $\hat{f}_n^{\lambda_n^{\alpha,\beta}}$, we need to know parameters $\alpha$ and $\beta$ that are not available in practice. Thus, it is important to construct a classifier, free from these parameters, which has the same behavior as $\hat{f}_n^{\lambda_n^{\alpha,\beta}}$, if the underlying distribution $\pi$ belongs to $\mathcal{Q}_{\alpha,\beta}$. Below we give such a construction.

Since the RKHS $H$ is given, the implementation of the $L1$-SVM classifier $\hat{f}_n^\lambda$ requires only knowledge of the regularization parameter $\lambda$. Thus, to provide an easily implemented procedure, using our aggregation method, it is natural to combine $L1$-SVM classifiers constructed for different values of $\lambda$ in a finite grid. We now define such a procedure.

We consider the $L1$-SVM classifiers $\hat{f}_m^\lambda$, defined in (3.3) for the subsample $D_m^1$, where $\lambda$ lies in the grid

$$\mathcal{G}(l) = \{\lambda_{l,k} = l^{-\phi_{l,k}} : \phi_{l,k} = 1/2 + k\Delta^{-1}, k = 0, \ldots, \lfloor 3\Delta/2 \rfloor\},$$

where we set $\Delta = l^{b_0}$ with some $b_0 > 0$. The subsample $D_l^2$ is used to aggregate these classifiers by the procedure (2.1), namely,

$$(3.5) \qquad \tilde{f}_n = \sum_{\lambda \in \mathcal{G}(l)} w_\lambda^{(l)} \hat{f}_m^\lambda,$$

where

$$w_\lambda^{(l)} = \frac{\exp(\sum_{i=m+1}^n Y_i \hat{f}_m^\lambda(X_i))}{\sum_{\lambda' \in \mathcal{G}(l)} \exp(\sum_{i=m+1}^n Y_i \hat{f}_m^{\lambda'}(X_i))} = \frac{\exp(-lA^{[l]}(\hat{f}_m^\lambda))}{\sum_{\lambda' \in \mathcal{G}(l)} \exp(-lA^{[l]}(\hat{f}_m^{\lambda'}))}$$

and $A^{[l]}(f) = (1/l)\sum_{i=m+1}^n (1 - Y_i f(X_i))_+$.

THEOREM 3.2. *Let $H$ be an RKHS with a continuous kernel on a compact metric space $\mathcal{X}$ satisfying* (A2) *with parameter $0 < p < 2$. Let $K$ be a compact subset of $(0,+\infty) \times (0,1]$. The classifier $\tilde{f}_n$, defined in* (3.5), *satisfies*

$$\sup_{\pi \in \mathcal{Q}_{\alpha,\beta}} \mathbb{E}[R(\tilde{f}_n) - R^*] \leq C n^{-4\beta(\alpha+1)/((2\alpha+p\alpha+4)(1+\beta))+\varepsilon}$$

for all $(\alpha, \beta) \in K$ and $\varepsilon > 0$, where $\mathcal{Q}_{\alpha,\beta}$ is the set of all probability measures on $\mathcal{X} \times \{-1, 1\}$ satisfying (MA2)$(\alpha)$ and (A2) with complexity parameter $(p, \beta)$ and $C > 0$ is a constant depending only on $\varepsilon, p, K, a$ and $b_0$.

3.2.2. *The case of the Gaussian RBF kernel.* In this subsection we apply our aggregation procedure to $L$1-SVM classifiers using the *Gaussian RBF kernel*. Let $\mathcal{X}$ be the closed unit ball of the space $\mathbb{R}^{d_0}$ endowed with the Euclidean norm $\|x\| = (\sum_{i=1}^{d_0} x_i^2)^{1/2}, \forall x = (x_1, \ldots, x_{d_0}) \in \mathbb{R}^{d_0}$. The Gaussian RBF kernel is defined as $K_\sigma(x, x') = \exp(-\sigma^2 \|x - x'\|^2)$ for $x, x' \in \mathcal{X}$, where $\sigma$ is a parameter and $\sigma^{-1}$ is called the *width* of the Gaussian kernel. The RKHS associated with $K_\sigma$ is denoted by $H_\sigma$.

Steinwart and Scovel [39] introduced the following assumption.

ASSUMPTION (GNA) (*Geometric noise assumption*). There exist $C_1 > 0$ and $\gamma > 0$ such that

$$\mathbb{E}\left[|2\eta(X) - 1| \exp\left(-\frac{\tau(X)^2}{t}\right)\right] \leq C_1 t^{\gamma d_0/2} \qquad \forall t > 0.$$

Here $\tau$ is a function on $\mathcal{X}$ with values in $\mathbb{R}$ which measures the distance between a given point $x$ and the decision boundary, namely,

$$\tau(x) = \begin{cases} d(x, G_0 \cup G_1), & \text{if } x \in G_{-1}, \\ d(x, G_0 \cup G_{-1}), & \text{if } x \in G_1, \\ 0, & \text{otherwise,} \end{cases}$$

for all $x \in \mathcal{X}$, where $G_0 = \{x \in \mathcal{X} : \eta(x) = 1/2\}$, $G_1 = \{x \in \mathcal{X} : \eta(x) > 1/2\}$ and $G_{-1} = \{x \in \mathcal{X} : \eta(x) < 1/2\}$. Here $d(x, A)$ denotes the Euclidean distance from a point $x$ to the set $A$. If $\pi$ satisfies Assumption (GNA) for a $\gamma > 0$, we say that $\pi$ has a *geometric noise exponent* $\gamma$.

The $L$1-SVM classifier associated to the Gaussian RBF kernel with width $\sigma^{-1}$ and regularization parameter $\lambda$ is defined by $\hat{f}_n^{(\sigma,\lambda)} = \text{sign}(\hat{F}_n^{(\sigma,\lambda)})$, where $\hat{F}_n^{(\sigma,\lambda)}$ is given by (3.3) with $H = H_\sigma$. Using the standard development related to SVM (cf. Schölkopf and Smola [37]), we may write $\hat{F}_n^{(\sigma,\lambda)}(x) = \sum_{i=1}^n \hat{C}_i K_\sigma(X_i, x), \forall x \in \mathcal{X}$, where $\hat{C}_1, \ldots, \hat{C}_n$ are solutions of the maximization problem

$$\max_{0 \leq 2\lambda C_i Y_i \leq n^{-1}} \left\{ 2 \sum_{i=1}^n C_i Y_i - \sum_{i,j=1}^n C_i C_j K_\sigma(X_i, X_j) \right\},$$

which can be obtained using standard quadratic programming software. According to Steinwart and Scovel [39], if the probability measure $\pi$ on $\mathcal{X} \times \{-1, 1\}$ satisfies the margin Assumption (MA2)$(\alpha)$ with margin parameter $0 \leq \alpha < +\infty$ and Assumption (GNA) with a geometric noise exponent $\gamma > 0$, the classifier $\hat{f}_n^{(\sigma_n^{\alpha,\gamma}, \lambda_n^{\alpha,\gamma})}$, where the regularization parameter and

width are defined by

$$\lambda_n^{\alpha,\gamma} = \begin{cases} n^{-(\gamma+1)/(2\gamma+1)}, & \text{if } \gamma \leq \dfrac{\alpha+2}{2\alpha}, \\ n^{-2(\gamma+1)(\alpha+1)/(2\gamma(\alpha+2)+3\alpha+4)}, & \text{otherwise} \end{cases}$$

and

$$\sigma_n^{\alpha,\gamma} = (\lambda_n^{\alpha,\gamma})^{-1/(\gamma+1)d_0},$$

satisfies

(3.6)
$$\mathbb{E}[R(\hat{f}_n^{(\sigma_n^{\alpha,\gamma},\lambda_n^{\alpha,\gamma})}) - R^*]$$
$$\leq C \begin{cases} n^{-\gamma/(2\gamma+1)+\varepsilon}, & \text{if } \gamma \leq \dfrac{\alpha+2}{2\alpha}, \\ n^{-2\gamma(\alpha+1)/(2\gamma(\alpha+2)+3\alpha+4)+\varepsilon}, & \text{otherwise}, \end{cases}$$

for all $\varepsilon > 0$, where $C > 0$ is a constant which depends only on $\alpha, \gamma$ and $\varepsilon$. We remark that fast rates are obtained only for $\gamma > (3\alpha + 4)/(2\alpha)$.

To construct the classifier $\hat{f}_n^{(\sigma_n^{\alpha,\gamma},\lambda_n^{\alpha,\gamma})}$, we need to know parameters $\alpha$ and $\gamma$, which are not available in practice. As in Section 3.2.1, we use our procedure to obtain a classifier which is adaptive to the margin and to the geometric noise parameters. Our aim is to provide an easily computed adaptive classifier. We propose the following method based on a grid for $(\sigma, \lambda)$. We consider the finite sets

$$\mathcal{M}(l) = \left\{ (\varphi_{l,p_1}, \psi_{l,p_2}) = \left( \frac{p_1}{2\Delta}, \frac{p_2}{\Delta} + \frac{1}{2} \right) : p_1 = 1, \ldots, 2\lfloor \Delta \rfloor; \right.$$
$$\left. p_2 = 1, \ldots, \lfloor \Delta/2 \rfloor \right\},$$

where we let $\Delta = l^{b_0}$ for some $b_0 > 0$, and

$$\mathcal{N}(l) = \{(\sigma_{l,\varphi}, \lambda_{l,\psi}) = (l^{\varphi/d_0}, l^{-\psi}) : (\varphi, \psi) \in \mathcal{M}(l)\}.$$

We construct the family of classifiers $(\hat{f}_m^{(\sigma,\lambda)} : (\sigma, \lambda) \in \mathcal{N}(l))$ using the observations of the subsample $D_m^1$ and we aggregate them by the procedure (2.1) using $D_l^2$, namely,

(3.7)
$$\tilde{f}_n = \sum_{(\sigma,\lambda) \in \mathcal{N}(l)} w_{\sigma,\lambda}^{(l)} \hat{f}_m^{(\sigma,\lambda)},$$

where

(3.8) $$w_{\sigma,\lambda}^{(l)} = \frac{\exp(\sum_{i=m+1}^n Y_i \hat{f}_m^{(\sigma,\lambda)}(X_i))}{\sum_{(\sigma',\lambda') \in \mathcal{N}(l)} \exp(\sum_{i=m+1}^n Y_i \hat{f}_m^{(\sigma',\lambda')}(X_i))} \qquad \forall (\sigma, \lambda) \in \mathcal{N}(l).$$

Denote by $\mathcal{R}_{\alpha,\gamma}$ the set of all probability measures on $\mathcal{X} \times \{-1,1\}$ satisfying both the margin Assumption (MA2)($\alpha$) with a margin parameter $\alpha > 0$ and Assumption (GNA) with a geometric noise exponent $\gamma > 0$. Define $\mathcal{U} = \{(\alpha,\gamma) \in (0,+\infty)^2 : \gamma > \frac{\alpha+2}{2\alpha}\}$ and $\mathcal{U}' = \{(\alpha,\gamma) \in (0,+\infty)^2 : \gamma \leq \frac{\alpha+2}{2\alpha}\}$.

THEOREM 3.3. *Let $K$ be a compact subset of $\mathcal{U}$ and $K'$ a compact subset of $\mathcal{U}'$. The aggregate $\tilde{f}_n$, defined in (3.7), satisfies*

$$\sup_{\pi \in \mathcal{R}_{\alpha,\gamma}} \mathbb{E}[R(\tilde{f}_n) - R^*] \leq C \begin{cases} n^{-\gamma/(2\gamma+1)+\varepsilon}, & \text{if } (\alpha,\gamma) \in K', \\ n^{-2\gamma(\alpha+1)/(2\gamma(\alpha+2)+3\alpha+4)+\varepsilon}, & \text{if } (\alpha,\gamma) \in K, \end{cases}$$

*for all $(\alpha,\gamma) \in K \cup K'$ and $\varepsilon > 0$, where $C > 0$ depends only on $\varepsilon, K, K', a$ and $b_0$.*

## 4. Proofs.

LEMMA 4.1. *For all positive $v, t$ and all $\kappa \geq 1$, $t + v \geq v^{(2\kappa-1)/2\kappa} t^{1/(2\kappa)}$.*

PROOF. Since log is concave, we have $\log(ab) = (1/x)\log(a^x) + (1/y) \times \log(b^y) \leq \log(a^x/x + b^y/y)$ for all positive numbers $a, b$ and $x, y$ such that $1/x + 1/y = 1$; thus $ab \leq a^x/x + b^y/y$. Lemma 4.1 follows by applying this relation with $a = t^{1/(2\kappa)}, x = 2\kappa$ and $b = v^{(2\kappa-1)/(2\kappa)}$. □

PROOF OF PROPOSITION 2.1. Observe that $(1-x)_+ = 1-x$ for $x \leq 1$. Since $Y_i \tilde{f}_n(X_i) \leq 1$ and $Y_i f_j(X_i) \leq 1$ for all $i = 1, \ldots, n$ and $j = 1, \ldots, M$, we have $A_n(\tilde{f}_n) = \sum_{j=1}^M w_j^{(n)} A_n(f_j)$. We have $A_n(f_j) = A_n(f_{j_0}) + \frac{1}{n}(\log(w_{j_0}^{(n)}) - \log(w_j^{(n)}))$, for any $j, j_0 = 1, \ldots, M$, where the weights $w_j^{(n)}$ are defined in (2.3) by

$$w_j^{(n)} = \frac{\exp(-nA_n(f_j))}{\sum_{k=1}^M \exp(-nA_n(f_k))},$$

and by multiplying the last equation by $w_j^{(n)}$ and summing over $j$, we get

(4.1) $$A_n(\tilde{f}_n) \leq \min_{j=1,\ldots,M} A_n(f_j) + \frac{\log M}{n}.$$

Indeed, we have $\log(w_{j_0}^{(n)}) \leq 0, \forall j_0 = 1, \ldots, M$, and $\sum_{j=1}^M w_j^{(n)} \log(\frac{w_j^{(n)}}{1/M}) = K(w|u) \geq 0$, where $K(w|u)$ denotes the Kullback–Leiber divergence between the weights $w = (w_j^{(n)})_{j=1,\ldots,M}$ and uniform weights $u = (1/M)_{j=1,\ldots,M}$. □

PROOF OF PROPOSITION 2.2. Denote $\gamma = (\log M)^{-1/4}$, $u = 2\gamma n^{-\kappa/(2\kappa-1)} \times \log^2 M$ and $W_n = (1-\gamma)(A(\tilde{f}_n) - A^*) - (A_n(\tilde{f}_n) - A_n(f^*))$. We have

$$\mathbb{E}[W_n] = \mathbb{E}[W_n(\mathbb{1}_{(W_n \leq u)} + \mathbb{1}_{(W_n > u)})]$$

$$\leq u + \mathbb{E}[W_n \mathbb{1}_{(W_n > u)}]$$

$$= u + u\mathbb{P}(W_n > u) + \int_u^{+\infty} \mathbb{P}(W_n > t)\,dt$$

$$\leq 2u + \int_u^{+\infty} \mathbb{P}(W_n > t)\,dt.$$

On the other hand, $(f_j)_{j=1,\ldots,M}$ are prediction rules, so we have $A(f_j) = 2R(f_j)$ and $A_n(f_j) = 2R_n(f_j)$ (recall that $A^* = 2R^*$). Moreover, we work in the linear part of the hinge-loss; thus

$$\mathbb{P}(W_n > t) = \mathbb{P}\left(\sum_{j=1}^M w_j((A(f_j) - A^*)(1 - \gamma) - (A_n(f_j) - A_n(f^*))) > t\right)$$

$$\leq \mathbb{P}\left(\max_{j=1,\ldots,M}((A(f_j) - A^*)(1 - \gamma) - (A_n(f_j) - A_n(f^*))) > t\right)$$

$$\leq \sum_{j=1}^M \mathbb{P}(Z_j > \gamma(R(f_j) - R^*) + t/2)$$

for all $t > u$, where $Z_j = R(f_j) - R^* - (R_n(f_j) - R_n(f^*))$ for all $j = 1, \ldots, M$ [recall that $R_n(f)$ is the empirical risk defined in (1.1)].

Let $j \in \{1, \ldots, M\}$. We can write $Z_j = (1/n)\sum_{i=1}^n (\mathbb{E}[\zeta_{i,j}] - \zeta_{i,j})$, where $\zeta_{i,j} = \mathbb{1}_{(Y_i f_j(X_i) \leq 0)} - \mathbb{1}_{(Y_i f^*(X_i) \leq 0)}$. We have $|\zeta_{i,j}| \leq 1$ and, under the margin assumption, we have $\mathbb{V}(\zeta_{i,j}) \leq \mathbb{E}(\zeta_{i,j}^2) = \mathbb{E}[|f_j(X) - f^*(X)|] \leq c(R(f_j) - R^*)^{1/\kappa}$, where $\mathbb{V}$ is the symbol of the variance. By applying Bernstein's inequality and Lemma 4.1 respectively, we get

$$\mathbb{P}[Z_j > \varepsilon] \leq \exp\left(-\frac{n\varepsilon^2}{2c(R(f_j) - R^*)^{1/\kappa} + 2\varepsilon/3}\right)$$

$$\leq \exp\left(-\frac{n\varepsilon^2}{4c(R(f_j) - R^*)^{1/\kappa}}\right) + \exp\left(-\frac{3n\varepsilon}{4}\right),$$

for all $\varepsilon > 0$. Denote $u_j = u/2 + \gamma(R(f_j) - R^*)$. After a standard calculation we get

$$\int_u^{+\infty} \mathbb{P}(Z_j > \gamma(R(f_j) - R^*) + t/2)\,dt = 2\int_{u_j}^{+\infty} \mathbb{P}(Z_j > \varepsilon)\,d\varepsilon \leq B_1 + B_2,$$

where

$$B_1 = \frac{4c(R(f_j) - R^*)^{1/\kappa}}{nu_j} \exp\left(-\frac{nu_j^2}{4c(R(f_j) - R^*)^{1/\kappa}}\right)$$

and

$$B_2 = \frac{8}{3n}\exp\left(-\frac{3nu_j}{4}\right).$$

Since $R(f_j) \geq R^*$, Lemma 4.1 yields $u_j \geq \gamma(R(f_j) - R^*)^{1/(2\kappa)}(\log M)^{(2\kappa-1)/\kappa} \times n^{-1/2}$. For any $a > 0$, the mapping $x \mapsto (ax)^{-1}\exp(-ax^2)$ is decreasing on $(0, +\infty)$; thus we have

$$B_1 \leq \frac{4c}{\gamma\sqrt{n}}(\log M)^{-(2\kappa-1)/\kappa}\exp\left(-\frac{\gamma^2}{4c}(\log(M))^{(4\kappa-2)/\kappa}\right).$$

The mapping $x \mapsto (2/a)\exp(-ax)$ is decreasing on $(0, +\infty)$ for any $a > 0$ and $u_j \geq \gamma(\log M)^2 n^{-\kappa/(2\kappa-1)}$; thus

$$B_2 \leq \frac{8}{3n}\exp\left(-\frac{3\gamma}{4}n^{(\kappa-1)/(2\kappa-1)}(\log M)^2\right).$$

Since $\gamma = (\log M)^{-1/4}$, we have $\mathbb{E}(W_n) \leq 4n^{-\kappa/(2\kappa-1)}(\log M)^{7/4} + T_1 + T_2$, where

$$T_1 = \frac{4Mc}{\sqrt{n}}(\log M)^{-(7\kappa-4)/(4\kappa)}\exp\left(-\frac{3}{4c}(\log M)^{(7\kappa-4)/(2\kappa)}\right)$$

and

$$T_2 = \frac{8M}{3n}\exp(-(3/4)n^{(\kappa-1)/(2\kappa-1)}(\log M)^{7/4}).$$

We have $T_2 \leq 6(\log M)^{7/4}/n$ for any integer $M \geq 1$. Moreover, $\kappa/(2\kappa - 1) \leq 1$ for all $1 \leq \kappa < +\infty$, so we get $T_2 \leq 6n^{-\kappa/(2\kappa-1)}(\log M)^{7/4}$ for any integers $n \geq 1$ and $M \geq 2$.

Let $B$ be a positive number. The inequality $T_1 \leq Bn^{-\kappa/(2\kappa-1)}(\log M)^{7/4}$ is equivalent to

$$2(2\kappa - 1)\left[\frac{3}{4c}(\log M)^{(7\kappa-4)/(2\kappa)} - \log M + \frac{7\kappa - 2}{2\kappa}\log(\log M)\right]$$
$$\geq \log((4c/B)^{2(2\kappa-1)}n).$$

Since we have $\frac{7\kappa-4}{2\kappa} \geq \frac{3}{2} > 1$ for all $1 \leq \kappa < +\infty$ and $M \geq an^b$ for some positive numbers $a$ and $b$, there exists a constant $B$ which depends only on $a, b$ and $c$ [for instance, $B = 4ca^{-1/(2b)}$ when $n$ satisfies $\log(an^b) \geq (b^2(8c/6)^2) \vee ((8c/3) \vee 1)^2$] such that $T_1 \leq Bn^{-\kappa/(2\kappa-1)}(\log M)^{7/4}$. □

PROOF OF THEOREM 2.3. Let $\gamma = (\log M)^{-1/4}$. Using (4.1), we have

$$\mathbb{E}[(A(\tilde{f}_n) - A^*)(1-\gamma)] - (A(f_{j_0}) - A^*)$$
$$= \mathbb{E}[(A(\tilde{f}_n) - A^*)(1-\gamma) - (A_n(\tilde{f}_n) - A_n(f^*))] + \mathbb{E}[A_n(\tilde{f}_n) - A_n(f_{j_0})]$$
$$\leq \mathbb{E}[(A(\tilde{f}_n) - A^*)(1-\gamma) - (A_n(\tilde{f}_n) - A_n(f^*))] + \frac{\log M}{n}.$$

For $W_n$ defined at the beginning of the proof of Proposition 2.2 and $f^*$ the Bayes rule, we have

$$(4.2) \quad (1-\gamma)(\mathbb{E}[A(\tilde{f}_n)] - A^*) \leq \min_{j=1,\ldots,M}(A(f_j) - A^*) + \mathbb{E}[W_n] + \frac{\log M}{n}.$$

According to Proposition 2.2, $\mathbb{E}[W_n] \leq C_0 n^{-\kappa/(2\kappa-1)}(\log M)^{7/4}$, where $C_0 > 0$ is given in Proposition 2.2. Using (4.2) and $(1-\gamma)^{-1} \leq 1 + 2\gamma$ for any $0 < \gamma < 1/2$, we get

$$\mathbb{E}[A(\tilde{f}_n) - A^*] \leq \left(1 + \frac{2}{\log^{1/4}(M)}\right)\left\{\min_{j=1,\ldots,M}(A(f_j) - A^*) + C\frac{\log^{7/4}(M)}{n^{\kappa/(2\kappa-1)}}\right\}.$$

We complete the proof by using inequality (1.3) and equality $2(R(f) - R^*) = A(f) - A^*$, which holds for any prediction rule $f$. □

PROOF OF THEOREM 2.4. Since the $\tilde{f}_k$'s take their values in $[-1,1]$ and $x \mapsto (1-x)_+$ is linear on $[-1,1]$, we obtain $A(\bar{f}_n) - A^* = \frac{1}{n}\sum_{k=1}^n (A(\tilde{f}_k) - A^*)$. Applying Theorem 2.3 to every $\tilde{f}_k$ for $k = 1,\ldots,n$, then taking the average of the $n$ oracle inequalities satisfied by the $\tilde{f}_k$ for $k = 1,\ldots,n$ and seeing that $(1/n)\sum_{k=1}^n k^{-\kappa/(2\kappa-1)} \leq \gamma(n,\kappa)$, we obtain

$$\mathbb{E}[A(\bar{f}_n) - A^*]$$
$$\leq \left(1 + \frac{2}{\log^{1/4}(M)}\right)\left\{\min_{j=1,\ldots,M}(A(f_j) - A^*) + C\gamma(n,\kappa)\log^{7/4}(M)\right\}.$$

We complete the proof by the same argument as at the end of the previous proof. □

PROOF OF THEOREM 3.1. Let $\rho_{\min} \leq \rho \leq \rho_{\max}$ and $\kappa \geq 1$. Let $\rho_{m,j_0} = \min(\rho_{m,j} : \rho_{m,j} \geq \rho)$. Since $N(m) \geq A'_0 m^{b'} \geq Cl^{b'}$, where $C > 0$, using the oracle inequality, stated in Theorem 2.3, we have, for $\pi$ satisfying (MA1)($\kappa$),

$$\mathbb{E}[R(\tilde{f}_n) - R^*|D_m^1]$$
$$\leq \left(1 + \frac{2}{\log^{1/4} N(m)}\right)\left\{2\min_{j=1,\ldots,N(m)}(R(\hat{f}_m^j) - R^*) + C\frac{\log^{7/4} N(m)}{l^{\kappa/(2\kappa-1)}}\right\},$$

where $C$ is a positive number depending only on $b', a, A'_0$ and $c$. Taking the expectation with respect to the subsample $D_m^1$, we have

$$\mathbb{E}[R(\tilde{f}_n) - R^*]$$
$$\leq \left(1 + \frac{2}{\log^{-1/4} N(m)}\right)\left\{2\mathbb{E}[R(\hat{f}_m^{j_0}) - R^*] + C\frac{\log^{7/4} N(m)}{l^{\kappa/(2\kappa-1)}}\right\}.$$

It follows from Tsybakov [42] that the excess risk of $\hat{f}_m^{j_0}$ satisfies

$$\sup_{\pi \in \mathcal{P}_{\kappa,\rho_{j_0}}} \mathbb{E}[R(\hat{f}_m^{j_0}) - R^*] \leq Cm^{-\kappa/(2\kappa+\rho_{j_0}-1)},$$

where $C$ is a positive number depending only on $A, c, \kappa, \rho_{\min}$ and $\rho_{\max}$ (note that $C$ does not depend on $\rho_{j_0}$).

Moreover, we have $m \geq n(1 - a/\log 3 - 1/3)$, $N(m) \leq A_0 m^b \leq A_0 n^b$ and $l \geq an/\log n$, so that there exists a constant $C$ depending only on $a, A_0, A_0'$, $b, b', \kappa, \rho_{\min}$ and $\rho_{\max}$ such that

(4.3) $$\sup_{\pi \in \mathcal{P}_{\kappa,\rho_{j_0}}} \mathbb{E}[R(\tilde{f}_n) - R^*] \leq C\{n^{-\kappa/(2\kappa+\rho_{j_0}-1)} + n^{-\kappa/(2\kappa-1)}(\log n)^{11/4}\}.$$

Since $\rho_{j_0} \leq \rho + N(m)^{-1} \leq \rho + (A_0')^{-1}[n(1 - a/\log 3 - 1/3)]^{-b'}$, there exists a constant $C$ depending only on $a, A_0', b', \kappa, \rho_{\min}$ and $\rho_{\max}$ such that, for all integers $n$, $n^{-\kappa/(2\kappa+\rho_{j_0}-1)} \leq Cn^{-\kappa/(2\kappa+\rho-1)}$. Theorem 2.4 follows directly from (4.3), seeing that $\rho \geq \rho_{\min} > 0$ and $\mathcal{P}_{\kappa,\rho} \subseteq \mathcal{P}_{\kappa,\rho_{j_0}}$ since $\rho_{j_0} \geq \rho$. □

PROOF OF THEOREM 3.2. Define $0 < \alpha_{min} < \alpha_{\max} < +\infty$ and $0 < \beta_{\min} < 1$ such that $K \subset [\alpha_{\min}, \alpha_{\max}] \times [\beta_{\min}, 1]$. Let $(\alpha_0, \beta_0) \in K$. We consider the function on $(0, +\infty) \times (0, 1]$ with values in $(1/2, 2)$, $\phi(\alpha, \beta) = 4(\alpha+1)/((2\alpha+p\alpha+4)(1+\beta))$. We take $k_0 \in \{0, \ldots, \lfloor 3\Delta/2 \rfloor - 1\}$ such that

$$\phi_{l,k_0} = 1/2 + k_0 \Delta^{-1} \leq \phi(\alpha_0, \beta_0) < 1/2 + (k_0+1)\Delta^{-1}.$$

For $n$ greater than a constant depending only on $K, p, b_0$ and $a$, there exists $\bar{\alpha}_0 \in [\alpha_{\min}/2, \alpha_{\max}]$ such that $\phi(\bar{\alpha}_0, \beta_0) = \phi_{l,k_0}$. Since $\alpha \mapsto \phi(\alpha, \beta_0)$ increases on $\mathbb{R}^+$, we have $\bar{\alpha}_0 \leq \alpha_0$. Moreover, we have $|\phi(\alpha_1, \beta_0) - \phi(\alpha_2, \beta_0)| \geq A|\alpha_1 - \alpha_2|, \forall \alpha_1, \alpha_2 \in [\alpha_{\min}/2, \alpha_{\max}]$, where $A > 0$ depends only on $p$ and $\alpha_{\max}$. Thus, $|\bar{\alpha}_0 - \alpha_0| \leq (A\Delta)^{-1}$. Since $\bar{\alpha}_0 \leq \alpha_0$, we have $\mathcal{Q}_{\alpha_0,\beta_0} \subseteq \mathcal{Q}_{\bar{\alpha}_0,\beta_0}$, so

$$\sup_{\pi \in \mathcal{Q}_{\alpha_0,\beta_0}} \mathbb{E}[R(\tilde{f}_n) - R^*] \leq \sup_{\pi \in \mathcal{Q}_{\bar{\alpha}_0,\beta_0}} \mathbb{E}[R(\tilde{f}_n) - R^*].$$

Since $\lceil 3\Delta/2 \rceil \geq (3/2)l^{b_0}$, for $\pi$ satisfying the margin Assumption (MA2)($\bar{\alpha}_0$), Theorem 2.3 leads to

$\mathbb{E}[R(\tilde{f}_n) - R^*|D_m^1]$

$$\leq \left(1 + \frac{2}{\log^{1/4}(\lceil 3\Delta/2 \rceil)}\right)\left\{2 \min_{\lambda \in \mathcal{G}(l)} (R(\hat{f}_m^\lambda) - R^*) + C_0 \frac{\log^{7/4}(\lceil 3\Delta/2 \rceil)}{l^{(\bar{\alpha}_0+1)/(\bar{\alpha}_0+2)}}\right\},$$

for all integers $n \geq 1$, where $C_0 > 0$ depends only on $K, a$ and $b_0$. Therefore, taking the expectation w.r.t. the subsample $D_m^1$, we get

$$\mathbb{E}[R(\tilde{f}_n) - R^*] \leq C_1(\mathbb{E}[R(\hat{f}_m^{\lambda_{l,k_0}}) - R^*] + l^{(\bar{\alpha}_0+1)/(\bar{\alpha}_0+2)} \log^{7/4}(n)),$$

where $\lambda_{l,k_0} = l^{-\phi_{l,k_0}}$ and $C_1 > 0$ depends only on $K, a$ and $b_0$.

Set $\Gamma:(0,+\infty)\times(0,1]\mapsto\mathbb{R}^+$ defined by $\Gamma(\alpha,\beta)=\beta\phi(\alpha,\beta), \forall(\alpha,\beta)\in(0,+\infty)\times(0,1]$. According to Steinwart and Scovel [38], if $\pi\in\mathcal{Q}_{\bar{\alpha}_0,\beta_0}$, then for all $\varepsilon>0$, there exists $C>0$, a constant depending only on $K,p$ and $\varepsilon$, such that

$$\mathbb{E}[R(\hat{f}_m^{\lambda_{l,k_0}})-R^*]\leq Cm^{-\Gamma(\bar{\alpha}_0,\beta_0)+\varepsilon}.$$

We remark that $C$ does not depend on $\bar{\alpha}_0$ and $\beta_0$ since $(\bar{\alpha}_0,\beta_0)\in[\alpha_{\min}/2,\alpha_{\max}]\times[\beta_{\min},1]$ and that the constant multiplying the rate of convergence, stated in Steinwart and Scovel [38], is uniformly bounded over $(\alpha,\beta)$ belonging to a compact subset of $(0,+\infty)\times(0,1]$.

Let $\varepsilon>0$. Assume that $\pi\in\mathcal{Q}_{\alpha_0,\beta_0}$. We have $n(1-a/\log 3-1/3)\leq m\leq n$, $l\geq an/\log n$ and $\Gamma(\bar{\alpha}_0,\beta_0)\leq(\bar{\alpha}_0+1)/(\bar{\alpha}_0+2)\leq 1$. Therefore, there exist $C_2,C_2'>0$ depending only on $a,b_0,K,p$ and $\varepsilon$ such that, for any $n$ greater than a constant depending only on $\beta_{\min},a$ and $b_0$,

$$\mathbb{E}[R(\tilde{f}_n)-R^*]\leq C_2(n^{-\Gamma(\bar{\alpha}_0,\beta_0)+\varepsilon}+n^{-(\bar{\alpha}_0+1)/(\bar{\alpha}_0+2)}(\log n)^{11/4})$$
$$\leq C_2'n^{-\Gamma(\bar{\alpha}_0,\beta_0)+\varepsilon}.$$

Moreover, $\Gamma$ satisfies $|\Gamma(\bar{\alpha}_0,\beta_0)-\Gamma(\alpha_0,\beta_0)|\leq B\Delta^{-1}$, where $B$ depends only on $p$ and $\alpha_{min}$, and $(n^{B\Delta^{-1}})_{n\in\mathbb{N}}$ is upper bounded. This completes the proof. □

PROOF OF THEOREM 3.3. Let $(\alpha_0,\gamma_0)\in K\cup K'$. First assume that $(\alpha_0,\gamma_0)$ belongs to $K\subset\mathcal{U}$. We consider the set

$$\mathcal{S}=\{(\varphi,\psi)\in(0,1/2)\times(1/2,1):2-2\psi-\varphi>0\}.$$

Each point of $\mathcal{S}$ is associated with a margin parameter (3.2) and with a geometric noise exponent by the following functions on $\mathcal{S}$ with values in $(0,+\infty)$:

$$\bar{\alpha}(\varphi,\psi)=\frac{4\psi-2}{2-2\psi-\varphi}\quad\text{and}\quad\bar{\gamma}(\varphi,\psi)=\frac{\psi}{\varphi}-1.$$

We take $(\varphi,\psi)\in\mathcal{S}\cap\mathcal{M}(l)$ such that $\bar{\alpha}(\varphi,\psi)\leq\alpha_0$, $\bar{\gamma}(\varphi,\psi)\leq\gamma_0$, $\bar{\alpha}(\varphi,\psi)$ is close enough to $\alpha_0$, $\bar{\gamma}(\varphi,\psi)$ is close enough to $\gamma_0$ and $\bar{\gamma}(\varphi,\psi)>(\bar{\alpha}(\varphi,\psi)+2)/(2\bar{\alpha}(\varphi,\psi))$. Since $\gamma_0>(\alpha_0+2)/(2\alpha_0)$, there exists a solution $(\varphi_0,\psi_0)\in\mathcal{S}$ of the system of equations

(4.4) $$\begin{cases}\bar{\alpha}(\varphi,\psi)=\alpha_0,\\ \bar{\gamma}(\varphi,\psi)=\gamma_0.\end{cases}$$

For all integers $n$ greater than a constant depending only on $K,a$ and $b_0$, there exists $(p_{1,0},p_{2,0})\in\{1,\ldots,2\lfloor\Delta\rfloor\}\times\{2,\ldots,\lfloor\Delta/2\rfloor\}$ defined by

$$\varphi_{l,p_{1,0}}=\min(\varphi_{l,p}:\varphi_{l,p}\geq\varphi_0)\quad\text{and}\quad\psi_{l,p_{2,0}}=\max(\psi_{l,p_2}:\psi_{l,p_2}\leq\psi_0)-\Delta^{-1}.$$

We have $2-2\psi_{l,p_{2,0}} - \varphi_{l,p_{1,0}} > 0$. Therefore, $(\varphi_{l,p_{1,0}}, \psi_{l,p_{2,0}}) \in \mathcal{S} \cap \mathcal{M}(l)$. Define $\bar{\alpha}_0 = \bar{\alpha}(\varphi_{l,p_{1,0}}, \psi_{l,p_{2,0}})$ and $\bar{\gamma}_0 = \bar{\gamma}(\varphi_{l,p_{1,0}}, \psi_{l,p_{2,0}})$. Since $(\varphi_0, \psi_0)$ satisfies (4.4), we have

$$\psi_{l,p_{2,0}} + \frac{1}{\Delta} \leq \psi_0 = \frac{-\alpha_0}{2\alpha_0+4}\varphi_0 + \frac{1+\alpha_0}{2+\alpha_0} \leq \frac{-\alpha_0}{2\alpha_0+4}\left(\varphi_{l,p_{1,0}} - \frac{1}{2\Delta}\right) + \frac{1+\alpha_0}{2+\alpha_0}$$

and $(\alpha_0/(2\alpha_0+4))(2\Delta)^{-1} \leq \Delta^{-1}$; thus

$$\psi_{l,p_{2,0}} \leq -\frac{\alpha_0}{2\alpha_0+4}\varphi_{l,p_{1,0}} + \frac{1+\alpha_0}{2+\alpha_0} \qquad \text{so } \bar{\alpha}_0 \leq \alpha_0.$$

With a similar argument, we have $\psi_{l,p_{2,0}} \leq (\alpha_0+1)\varphi_{l,p_{1,0}}$, that is, $\bar{\gamma}_0 \leq \gamma_0$. Now we show that $\bar{\gamma}_0 > (\bar{\alpha}_0+2)/(2\bar{\alpha}_0)$. Since $(\alpha_0, \gamma_0)$ belongs to a compact, $(\varphi_0, \psi_0)$ and $(\varphi_{l,p_{1,0}}, \psi_{l,p_{2,0}})$ belong to a compact subset of $(0,1/2) \times (1/2,1)$ for $n$ greater than a constant depending only on $K, a, b_0$. Thus, there exists $A > 0$, depending only on $K$, such that, for $n$ large enough, we have

$$|\alpha_0 - \bar{\alpha}_0| \leq A\Delta^{-1} \quad \text{and} \quad |\gamma_0 - \bar{\gamma}_0| \leq A\Delta^{-1}.$$

Denote $d_K = d(\partial \mathcal{U}, K)$, where $\partial \mathcal{U}$ is the boundary of $\mathcal{U}$ and $d(A,B)$ denotes the Euclidean distance between sets $A$ and $B$. We have $d_K > 0$ since $K$ is a compact, $\partial \mathcal{U}$ is closed and $K \cap \partial \mathcal{U} = \varnothing$. Set $0 < \alpha_{\min} < \alpha_{\max} < +\infty$ and $0 < \gamma_{\min} < \gamma_{\max} < +\infty$ such that $K \subset [\alpha_{\min}, \alpha_{\max}] \times [\gamma_{\min}, \gamma_{\max}]$. Define $\mathcal{U}_\mu = \{(\alpha, \gamma) \in (0, +\infty)^2 : \alpha \geq 2\mu \text{ and } \gamma > (\alpha - \mu + 2)/(2(\alpha - \mu))\}$ for $\mu = \min(\alpha_{\min}/2, d_K)$. We have $K \subset \mathcal{U}_\mu$, so $\gamma_0 > (\alpha_0 - \mu + 2)/(2(\alpha_0 - \mu))$. Since $\alpha \mapsto (\alpha+2)/(2\alpha)$ is decreasing, $\bar{\gamma}_0 > \gamma_0 - A\Delta^{-1}$ and $\alpha_0 \leq \bar{\alpha}_0 + A\Delta^{-1}$, we have $\bar{\gamma}_0 > \bar{\beta}(\bar{\alpha}_0) - A\Delta^{-1}$, where $\bar{\beta}$ is a positive function on $(0, 2\alpha_{\max}]$ defined by $\bar{\beta}(\alpha) = (\alpha - (\mu - A\Delta^{-1}) + 2)/(2(\alpha - (\mu - A\Delta^{-1})))$. We have $|\bar{\beta}(\alpha_1) - \bar{\beta}(\alpha_2)| \geq (2\alpha_{\max})^{-2}|\alpha_1 - \alpha_2|$ for all $\alpha_1, \alpha_2 \in (0, 2\alpha_{\max}]$. Therefore, $\bar{\beta}(\bar{\alpha}_0) - A\Delta^{-1} \geq \bar{\beta}(\bar{\alpha}_0 + 4A\alpha_{\max}^2 \Delta^{-1})$. Thus, for $n$ greater than a constant depending only on $K, a$ and $b_0$, we have $\bar{\gamma}_0 > (\bar{\alpha}_0 + 2)/(2\bar{\alpha}_0)$.

Since $\bar{\alpha}_0 \leq \alpha_0$ and $\bar{\gamma}_0 \leq \gamma_0$, we have $\mathcal{R}_{\alpha_0,\gamma_0} \subset \mathcal{R}_{\bar{\alpha}_0,\bar{\gamma}_0}$ and

$$\sup_{\pi \in \mathcal{R}_{\alpha_0,\gamma_0}} \mathbb{E}[R(\tilde{f}_n) - R^*] \leq \sup_{\pi \in \mathcal{R}_{\bar{\alpha}_0,\bar{\gamma}_0}} \mathbb{E}[R(\tilde{f}_n) - R^*].$$

If $\pi$ satisfies (MA2)($\bar{\alpha}_0$), then we get from Theorem 2.3

$$\mathbb{E}[R(\tilde{f}_n) - R^* | D_m^1]$$

(4.5)
$$\leq \left(1 + \frac{2}{\log^{1/4} M(l)}\right)$$

$$\times \left\{2 \min_{(\sigma,\lambda) \in \mathcal{N}(l)} (R(\hat{f}_m^{(\sigma,\lambda)}) - R^*) + C_2 \frac{\log^{7/4}(M(l))}{l^{(\bar{\alpha}_0+1)/(\bar{\alpha}_0+2)}}\right\},$$

for all integers $n \geq 1$, where $C_2 > 0$ depends only on $K, a$ and $b_0$ and $M(l)$ is the cardinality of $\mathcal{N}(m)$. We remark that $M(l) \geq l^{2b_0}/2$, so we can apply Theorem 2.3.

Let $\varepsilon > 0$. Since $M(l) \leq n^{2b_0}$ and $\bar{\gamma}_0 > (\bar{\alpha}_0 + 2)/(2\bar{\alpha}_0)$, taking expectations in (4.5) and using the result (3.6) of Steinwart and Scovel [39], for $\sigma = \sigma_{l,\varphi_{l,p_{1,0}}}$ and $\lambda = \lambda_{l,\psi_{l,p_{2,0}}}$, we obtain

$$\sup_{\pi \in \mathcal{R}_{\bar{\alpha}_0, \bar{\gamma}_0}} \mathbb{E}[R(\tilde{f}_n) - R^*] \leq C(m^{-\Theta(\bar{\alpha}_0, \bar{\gamma}_0) + \varepsilon} + l^{-(\bar{\alpha}_0+1)/(\bar{\alpha}_0+2)} \log^{7/4}(n)),$$

where $\Theta : \mathcal{U} \mapsto \mathbb{R}$ is defined for all $(\alpha, \gamma) \in \mathcal{U}$ by $\Theta(\alpha, \gamma) = (2\gamma(\alpha+1))/(2\gamma(\alpha+2) + 3\alpha + 4)$ and $C > 0$ depends only on $a, b_0, K$ and $\varepsilon$. We remark that the constant before the rate of convergence in (3.6) is uniformly bounded on every compact of $\mathcal{U}$. We have $\Theta(\bar{\alpha}_0, \bar{\gamma}_0) \leq \Theta(\alpha_0, \gamma_0) \leq \Theta(\bar{\alpha}_0, \bar{\gamma}_0) + 2A\Delta^{-1}$, $m \geq n(1 - a/\log 3 - 1/3)$ and $(m^{2A\Delta^{-1}})_{n \in \mathbb{N}}$ is upper bounded, so there exists $C_1 > 0$ depending only on $K, a, b_0$ such that $m^{-\Theta(\bar{\alpha}_0, \bar{\gamma}_0)} \leq C_1 n^{-\Theta(\alpha_0, \gamma_0)}, \forall n \geq 1$.

A similar argument as at the end of the proof of Theorem 3.2 and the fact that $\Theta(\alpha, \gamma) < (\alpha+1)/(\alpha+2)$ for all $(\alpha, \gamma) \in \mathcal{U}$ lead to the result of the first part of Theorem 3.3.

Let now $(\alpha_0, \gamma_0) \in K'$. Let $\alpha'_{\max} > 0$ be such that $\forall (\alpha, \gamma) \in K', \alpha \leq \alpha'_{\max}$. Take $p_{1,0} \in \{1, \ldots, 2\lfloor \Delta \rfloor\}$ such that $\varphi_{l, p_{1,0}} = \min(\varphi_{l,p} : \varphi_{l,p} \geq (2\gamma_0 + 1)^{-1}$ and $p \in 4\mathbb{N})$, where $4\mathbb{N}$ is the set of all integer multiples of 4. For large values of $n$, $p_{1,0}$ exists and $p_{1,0} \in 4\mathbb{N}$. Denoting $\bar{\gamma}_0 \in (0, +\infty)$ such that $\varphi_{l,p_{1,0}} = (2\bar{\gamma}_0 + 1)^{-1}$, we have $\bar{\gamma}_0 \leq \gamma_0$; thus $\mathcal{R}_{\alpha_0, \gamma_0} \subseteq \mathcal{R}_{\alpha_0, \bar{\gamma}_0}$ and

$$\sup_{\pi \in \mathcal{R}_{\alpha_0, \gamma_0}} \mathbb{E}[R(\tilde{f}_n) - R^*] \leq \sup_{\pi \in \mathcal{R}_{\alpha_0, \bar{\gamma}_0}} \mathbb{E}[R(\tilde{f}_n) - R^*].$$

If $\pi$ satisfies the margin assumption (3.2) with the margin parameter $\alpha_0$, then, using Theorem 2.3, we obtain, for any integer $n \geq 1$,

$$\mathbb{E}[R(\tilde{f}_n) - R^* | D_m^1]$$

(4.6)
$$\leq \left(1 + \frac{2}{\log^{1/4}(M(l))}\right)$$

$$\times \left\{ 2 \min_{(\sigma, \lambda) \in \mathcal{N}(l)} (R(\hat{f}_m^{(\sigma, \lambda)}) - R^*) + C_0 \frac{\log^{7/4} M(l)}{l^{(\alpha_0+1)/(\alpha_0+2)}} \right\},$$

where $C > 0$ appears in Proposition 2.2 and $M(l)$ is the cardinality of $\mathcal{N}(l)$.

Let $\varepsilon > 0$ and $p_{2,0} \in \{1, \ldots, \lfloor \Delta/2 \rfloor\}$ be defined by $p_{2,0} = p_{1,0}/4$ (note that $p_{1,0} \in 4\mathbb{N}$). We have

$$\sigma_{l, \varphi_{l, p_{1,0}}} = (\lambda_{l, \psi_{l, p_{2,0}}})^{-1/(d_0(\bar{\gamma}_0+1))}.$$

Since $\bar{\gamma}_0 \leq (\alpha_0 + 2)/(2\alpha_0)$, using (3.6) of Steinwart and Scovel [39], we have, for $\sigma = \sigma_{l, \varphi_{l, p_{1,0}}}$ and $\lambda = \lambda_{l, \psi_{l, p_{2,0}}}$,

$$\mathbb{E}[R(\hat{f}_m^{(\sigma_0, \lambda_0)}) - R^*] \leq Cm^{-\bar{\Gamma}(\bar{\gamma}_0) + \varepsilon},$$

where $\bar{\Gamma}:(0,+\infty)\longmapsto\mathbb{R}$ is the function defined by $\bar{\Gamma}(\gamma)=\gamma/(2\gamma+1)$ for all $\gamma\in(0,+\infty)$ and $C>0$ depends only on $a,b_0,K'$ and $\varepsilon$. We remark that, as in the first part of the proof, we can uniformly bound the constant before the rate of convergence in (3.6) on every compact subset of $\mathcal{U}'$. Since $M(l)\leq n^{2b_0}$, taking the expectation in (4.6), we find

$$\sup_{\pi\in\mathcal{R}_{\alpha_0,\bar{\gamma}_0}}\mathbb{E}[R(\tilde{f}_n)-R^*]\leq C(m^{-\Gamma(\bar{\gamma}_0)+\varepsilon}+l^{-(\alpha_0+1)/(\alpha_0+2)}\log^{7/4}(n)),$$

where $C>0$ depends only on $a,b_0,K'$ and $\varepsilon$. Moreover, $|\gamma_0-\bar{\gamma}_0|\leq 2(2\alpha'_{\max}+1)^2\Delta^{-1}$, so $|\bar{\Gamma}(\bar{\gamma}_0)-\bar{\Gamma}(\gamma_0)|\leq 2(2\alpha_{\max}+1)\Delta^{-1}$. To achieve the proof, we use the same argument as for the first part of the proof. $\square$

**Acknowledgments.** I want to thank my advisor, Alexandre Tsybakov, for spending many hours to help me with this work, giving me many ideas, much advice and many interesting papers.

LABORATOIRE DE PROBABILITÉS ET MODÈLES
ALÉATOIRES (UMR CNRS 7599)
UNIVERSITÉ PARIS VI
4, PL. JUSSIEU, BOÎTE COURRIER 188
75252 PARIS CEDEX 05
FRANCE
E-MAIL: [lecue@ccr.jussieu.fr](lecue@ccr.jussieu.fr)